\documentclass[twoside,12pt,a4paper]{amsart} 
\usepackage{amsmath,amsfonts,amsthm,eurosym}
\usepackage[colorlinks=true,linkcolor=blue,citecolor=blue,filecolor=blue,urlcolor=blue,pageanchor=true,plainpages=false,linktocpage]{hyperref}
\usepackage[english]{babel}

\numberwithin{equation}{section}

\newcommand{\R}{\mathbb{R}}

\renewcommand{\epsilon}{\varepsilon}
\renewcommand{\phi}{\varphi}
\renewcommand{\rho}{\varrho}
\renewcommand{\theta}{\vartheta}

\begin{document}
\baselineskip3.6ex


\title{}

\author{Alexander Pigazzini$^{(*)}$}

\date{}

%


%

%

\makeatletter
\def\@makefnmark{}
\makeatother
\newcommand{\myfootnote}[2]{\footnote{\textbf{#1}: #2}}
 \footnote {$^{(*)}$ IT-IMPRESA s.r.l., Monza, MB 20900, ITALY - email:  pigazzinialexander18@gmail.com}
{\bfseries\centerline{On the (2+2)-Einstein Warped Product Manifolds with $f$-curvature-Base}}
\\
\\
\centerline{Alexander Pigazzini $^{(*)}$}
\\
\\
{\bfseries \centerline {Abstract}} 
\\ 
\\
We study the $(2+2)$-Einstein warped product manifolds, where the scalar 
\\
curvature of the Base is a multiple of the warping function, and we called this condition (inside a warped product manifold) {\textit {$f$-curvature-Base}} ($R_{f_B}$).The aim of this paper is to check if there are Base-manifolds with non-flat metrics that satisfy this condition, and this was done in cases where $M$ and Fiber-manifold are not both non-Ricci-flat. As a results of our cases we find that the "{\textit {$f$-curvature-Base}}" is equivalent to requesting a flat metric on the Base-manifold.
\\
\\
{\textit {Keywords}}:  {\textit Special Einstein warped product manifolds, multiple of the warping 
\\
function, $f$-curvature-Base, $2$-dimensional Base, $2$-dimensional Fiber, $(2+2)$ signature.}
\\
\\
{\textit {2010 Mathematics Subject Classification}}: 53C25, 53C21
\\
\\
{\bfseries \centerline { Introduction and Preliminaries}}
\\
\\
The warped product manifolds are very interesting for Riemannian Geometry and also for General Relativity, as a model of space-time, and in the same that way the Einstein warped product manifolds are object of study for mathematics and physics.
In the last few years several authors, have dealt with the study of the Einstein warped product manifolds in which the Base-manifold was a $2$-dimensional manifold (see for example [2], [7] e [9]). In [2] the author shows that the nonconstant warping functions $f$, for 2-dimensional Base-manifold, exists only if the metrics of the Base-manifold is of the form: $g_B=dt^2+f'(t)^2du^2$. The aim of this paper is to study the Einstein warped product manifolds $(M, \bar{g})$ with constrained scalar curvature of the Base-manifold $(B,g)$, namely $R_B$ is a multiply of the warping function $f$. We try to see if this condition admits $(B,g)$ with non-trivial metrics which must be of the form considered in [2]. We have called this condition, inside the warped product manifold, {\textit {$f$-curvature-Base}} ($R_{f_B}$) just to indicate that the curvature of $(B,g)$ is bounded to $f$. We have studied this condition in the following cases: 
\begin{itemize}
\item [\ -]or only $(M, \bar{g})$ is Ricci-flat (i.e. $\bar{Ric}=\lambda \bar{g}$, with $\lambda=0$, case $(1b)$),
\item [\ -]or only Fiber-manifold ($F, \ddot{g}$) is Ricci-flat (i.e. $\ddot{Ric}=\mu \ddot{g}$, with $\mu=0$, case $(2b)$), 
\item [\ -]or both are Ricci-flat at the same time (case $(1a)$=case $(2a)$).
\end{itemize} 
Note that the condition {\textit {"$f$-curvature-Base"}} means to find a metric for the Base-manifold (we are interested in the non-flat case) that also allows the solution for a particular nonlinear elliptic pde of the form: $\Delta f +af^2=0$ or $\Delta f +af^2+bf=0$ (our equations are well known pde, for example in [14] the author studied the non-homogeneous case $\Delta f +af^2=v$, see [5] for a basic pde theory).
\\
This paper addresses the $(2+2)$-Einstein warped product manifolds, a type of signatures that, over the years, has been the subject of many studies both from mathematical point of view ([1], [15]) and from the physical point of view ([3], [4], [6], [11] and [12]).
\\
We have found that there are no Base-manifolds with non-trivial metrics for the cases examined.
In our analysis we exclude a priori the case where the warping function $f$ is null ($f=0$) as a trivial case. In fact it will not be considered in the initial hypotheses.
\\
\\
{\bfseries Definition 1:} From [10] (see also [2], [8] and [13]) a warped product manifold is Einstein if only if:
\\
\\
\numberwithin{equation}{section}
{(1)}
$\bar{Ric}=\lambda \bar{g} \Longleftrightarrow\begin{cases} 
 Ric- \frac{m}{f}\nabla^2 f= \lambda g  \\  \ddot{Ric}=\mu \ddot{g} \\ f \Delta f+(m-1) |\nabla f|^2 + \lambda f^2 =\mu
\end{cases}$
\\
\\
where $\lambda$ and $\mu$ are constants, $m$ is the dimension of $F$, $\nabla ^2f$, $ \Delta f$ and $\nabla f$ are, 
\\
respectively, the Hessian, the Laplacian and the gradient of $f$ for $g$, with $f:(B) \rightarrow (0, \infty)$ a smooth positive function.
\\
It follows from (1): 
\\
\numberwithin{equation}{section}
{(2)}
$R_Bf^2-mf \Delta f=n f^2 \lambda$ 
\\
where $n$ and $R_B$ are the dimension and the scalar curvature of $B$ respectively; and from third equation, considering $m>1$, we have:
\\
\\
\numberwithin{equation}{section}
{(3)}
$mf\Delta f+m(m-1)|\nabla f|^2+m\lambda f^2=m\mu$
\\
Now from (2) and (3) we obtain:
\\
\\
\numberwithin{equation}{section}
{(4)}
$|\nabla f|^2+[\frac{\lambda (m-n)+R_B}{m(m-1)}]f^2=\frac{\mu}{(m-1)}$
\\
{\bfseries Definition 2:} Let $(M, \bar{g})=(B,g)\times_f(F,\ddot{g})$ be an Einstein warped-product manifold with $\bar{g}=g+f^2 \ddot{g}$. We define the scalar curvature of the Base-manifold $(B,g)$ as {\textit {$f$-curvature-Base}} ($R_{f_B}$), if it is a multiple of the warping function $f$ (i.e. $R_{f_B}=cf$ for $c$ an arbitrary constant $\in \R$).
\\
\\
{\bfseries \centerline { Case 1:Einstein warped-product manifold Ricci-flat ($\lambda=0$)}}
\\
\\
{\bfseries Case $1a$:}{\bfseries } $\mu=0$.
\\
\\
{\bfseries Theorem 1.} Let $(B^2,g)$ be a smooth surface with $R_{f_B}$, and $(F^2,\ddot{g})$ a smooth Ricci-flat surface (i.e. $\ddot{Ric}=\mu \ddot{g}$, with $\mu=0$).
\\
If $(M^{2+2}, \bar{g})=(B^2,g)\times_f(F^2,\ddot{g})$ is an Einstein warped-product manifold 
\\
Ricci-flat (i.e. $\bar{Ric}=\lambda \bar{g}$, where $\lambda=0$), then $f$ can only be a constant function.
\\
\\
{\textit {Proof.}} In our case, we have $m=n=2$, $\mu=0$ and $R_B=R_{f_B}$. Then (2) and (3) 
\\
become:
\\
\\
\numberwithin{equation}{section}
{(5)}
$\Delta f-hf^2+\lambda f=0$ \; \; \; (with $h=c/2$.)
\\
\numberwithin{equation}{section}
{(6)} 
$f\Delta f+|\nabla f|^2+\lambda f^2=0$
\\
\\
Then (4) becomes:
\\
\numberwithin{equation}{section}
{(7)}
$|\nabla f|^2+hf^3=0$
\\
\\
If we consider $h=0$, it is easy to see that $f$ must be a constant. So we assume $h \neq 0$. By setting $u=-hf$, the above equations (5) and (6) (with $\lambda =0$) becomes:
\\
\numberwithin{equation}{section}
{(8)} 
$\Delta u+u^2=0$
\\
\numberwithin{equation}{section}
{(9)}
$u\Delta u+|\nabla u|^2=0$
\\
\\
Then (7) becomes:
\\
\\
\numberwithin{equation}{section}
{(10)}
$|\nabla u|^2-u^3=0$
\\
\\
Let $g$ be the metric on $B$ and assume that $u$ is a nonzero (and hence necessarily positive) solution, to the above system on a simply-connected open subset $B' \subset B$.
\\
The equation (10) implies that $\omega_1=u^{-3/2}du$ is a $1$-form with $g$-norm $1$ on $B'$ and hence $g$ can be written in the form $g=\omega_1^2+\omega_2^2$ on $B'$ for some $\omega_2$, which is also a unit $1$-form. Fix an orientation by requiring that $\omega_1 \wedge \omega_2$ be the $g$-area form on $B'$, then $\star du=u^{3/2} \omega_2$ and since $d(\star du)=\Delta u \; \omega_1 \wedge \omega_2$, it follows that:
\\
$\frac{3}{2}u^{1/2}du \wedge \omega_2+u^{3/2}d\omega_2=d(u^{3/2}\omega_2)=-u^2 \omega_1 \wedge \omega_2=-u^{1/2}du \wedge \omega_2$
\\
or $d\omega_2=-\frac{5}{2}u^{-1}du \wedge \omega_2$, which can be written as: $d(u^{5/2} \omega_2)=0$.
\\
\\
Since $B'$ is simply connected, it follows that there exists a function $v$ on $B'$ such that $u^{5/2} \omega_2=dv$. Consequently, the metric $g$ has the form:
\\
\\
\numberwithin{equation}{section}
{(11)}
$g=\omega_1^2+\omega_2^2=u^{-3}du^2+u^{-5}dv^2$.
\\
This metric can be expressed in polar form by setting $u=4r^{-2}$ and $v=32\theta$, in which case, it becomes:
\\
\\
\numberwithin{equation}{section}
{(12)}
$g=dr^2+r^{10}d\theta ^2$
\\
\\
This is a singular and incomplete metric at $r=0$, though it is complete at $r=\infty$ and its Gauss curvature is $K=-\frac{(r^5)''}{r^5}=-20r^{-2}=-5u<0$, the original $f$ takes the form: $f=\frac{-u}{h}=\frac{-4}{hr^2}$.
\\
\\
From the initial hypothesis, where we had set $R_B=R_{f_B}$, we have that $K$ is 
\\
incompatible, then there is no solutions if the background metric $g$ is not flat, 
\\
so the only solution that we have is for $g$ flat and this implies $R_{f_B}=0$, i.e. $h=0$ and $f=constant$.  \; \; \; \; \; \; \; $\Box$
\\
\\
{\bfseries Case $1b$:}{\bfseries } $\mu \neq 0$.
\\
\\
{\bfseries Theorem 2.} Let $(B^2,g)$ be a smooth surface with $R_{f_B}$, and $(F^2,\ddot{g})$ a smooth Einstein-surface (i.e. $\ddot{Ric}=\mu \ddot{g}$).
\\
If $(M^{2+2}, \bar{g})=(B^2,g)\times_f(F^2,\ddot{g})$ is an Einstein warped-product manifold 
\\
Ricci-flat (i.e. $\bar{Ric}=\lambda \bar{g}$, then $f$ can only be an affine function.
\\
\\
{\textit {Proof.}} From (1) and considering $\lambda=0$, the equations (5) and (7) become:
\\
\\
(13) $\Delta f - hf^2=0$
\\
(14) $|\nabla f|^2+hf^3 - \mu =0$
\\
Now, if we consider the flat-metric ($h=0$), it is easy to see that $f$ must be an affine function, while if $f$ is a constant, this implies $h=\mu=0$ (and $\mu=0$ is not the case in our interest); so we assume $h \neq 0$ with $f$ nonconstant and  setting $u=-hf$, for an open set where $u$ nonzero, then:
\\
(15) $\Delta u + u^2=0$
\\
(16) $|\nabla u|^2- u^3 - h^2\mu =0$.
\\
\\
Now, for semplicity we replace the constant $h^2 \mu$ with constant $A$.
\\
The equation (16) implies that $\omega_1=(u^3+A)^{-\frac{1}{2}}du$, and this implies that we have to assume $(u^3+A)$ nonzero, then $g = \omega_1^2 + \omega_2^2$ for some $\omega_2$. 
\\
Following the procedure used for {\textit {Theorem 1}}, we have: $\star du=(u^3+A)^{\frac{1}{2}} \omega_2$, and $d(\star du)=\Delta u \; \omega_1 \wedge \omega_2$. Now we have:
\\
\\
$\frac{3}{2}(u^3+A)^{-\frac{1}{2}} u^2 du \wedge \omega_2 + (u^3 + A)^{\frac{1}{2}} d\omega_2=$
\\
$=d[(u^3+A)^{\frac{1}{2}} \omega_2]=-u^2 \omega_1 \wedge \omega_2=-u^2(u^3+A)^{-\frac{1}{2}}du \wedge \omega_2=$
\\
$=-\frac{3}{2}(u^3+A)^{-1}u^2du \wedge \omega_2-u^2(u^3+A)^{-1}du \wedge \omega_2=d\omega_2$.
\\
\\
\\
Then $-\frac{5}{2}(u^3+A)^{-1}u^2 du \wedge \omega_2=d\omega_2$ and we have $\omega_2=(u^3+A)^{-\frac{5}{6}} dv$. Thus the 
\\
metric $g=(u^3+A)^{-1}du^2+ (u^3+A)^{-\frac{5}{3}}dv^2$ has a singularity in $u^3=-A$.
\\
\\
{\bfseries [Note 1:} The purpose of this analysis is to check if there is an incompatibility between the initial hypothesis and the Gaussian curvature built by the obtained metric and therefore does not concern the study, in depth, of singularities.
\\
So only some considerations will be made.{\bfseries ]}
\\
\\
Given that we had to assume that $u^3+A$ was nonzero, see $\omega_1$ obtained from (16), we also consider that $u^3+A$ is never zero, and to study what happens near this locus we need to insert the cases according to whether $A=0$ or not.  
\\
If $A=0$ then we have $\mu=0$ (and we have the same situation of the "case 1a"), but for the initial assumptions we must have $\mu \neq 0$, which implies that we can not consider the case $A=0$.
\\
Now, we assume that $A$ is not zero. By hypothesis, $u^3+A=|\nabla u|^2$ is always non-negative. Write $A=-a^3$ for some unique constant $a$. 
\\
\\
Then we have:
\\
\\
$|\nabla u|^2=(u-a)(u^2+ua+a^2)$ 
\\
and the metric is: $g=[(u-a)(u^2+ua+a^2)]^{-1}du^2+(u-a)(u^2+ua+a^2)^{-\frac{5}{3}}dv^2$.
\\
Since $a$ is non-zero, it follows that $(u^2+ua+a^2)$ is always positive, so if $u-a$ 
\\
vanishes along a gradient line (which must be a geodesic), then it must vanish exactly to order 2 along that gradient line.
\\
\\
There are two cases: First, $u-a$ vanishes at an isolated point, in which case our surface is rotationally symmetric about that point. This means that our coordinate system will break down 
\\
(just as polar coordinates do) or else, we rewrite $u=a+s^2$ for some function $s$ whose gradient along where it vanishes is nonzero.
\\
\\
At the end of these considerations, being that for the purposes of our work we have considered an open set for which $(u^3+A)$ is nonzero and the Gaussian curvature will be given by: $K=5u - 10u^4(u^3+A)^{-1}$.
\\
\\
Also in this case, it is easy to verify that for the initial hypothesis, where we had set $R_B=R_{f_B}$ (i.e. $K=-u$), we have that $K$ is incompatible with our analysis. In fact this implies that $f$ constant  (i.e. $h=\mu=0$). Then the only solution is a flat metric $g$, (i.e. $h=0$) which in this case implies that $f$ is an affine function. \; \; \; \; \; \; \; $\Box$
\\
\\
{\bfseries \centerline { Case 2: Flat Fiber Surface ($\mu=0$)}}
\\
\\
{\bfseries Case $2a$:}{\bfseries } We point out that if we consider $\lambda=0$ we are obviously in the same situation treated in the case $(1a)$, so we will not have to dwell further.
\\
\\
{\bfseries Case $2b$:}{\bfseries } $\lambda \neq 0$.
\\
{\bfseries [ Note 2}: If we set $h=0$ we force $\lambda=0$ and this is not the case in questions .{\bfseries ]}
\\
\\
{\bfseries Theorem 3.} Let $(B^2,g)$ be a smooth surface with $R_{f_B}$, and $(F^2,\ddot{g})$ a smooth Ricci-flat surface (i.e. $\ddot{Ric}=\mu \ddot{g}$, with $\mu=0$).
\\
If $(M^{2+2}, \bar{g})=(B^2,g)\times_f(F^2,\ddot{g})$ is an Einstein warped-product manifold (i.e. $\bar{Ric}=\lambda \bar{g}$), then $(M^{2+2}, \bar{g})$ exists if and only if $\lambda=0$.
\\
{\textit {Proof.}} The analysis is essentially the same as seen so far, so we assume $h \neq 0$ and set $u=-hf$. Now our equations become:
\\
\\
\numberwithin{equation}{section}
{(17)}
$\Delta u+u^2+\lambda u=0$
\\
\numberwithin{equation}{section}
{(18)} 
$u\Delta u+|\nabla u|^2+\lambda u^2=0$
\\
\\
Then from (17) and (18) we have:
\\
\\
\numberwithin{equation}{section}
{(19)} 
$|\nabla u|^2-u^3=0$
\\
\\
On the open set where $u \neq 0$ it is positive, so if we set $u=4/r^2$ we find that $\omega_1=dr$ and $\omega_2=r^5 e^{(\lambda /4)r^2}d\theta$ for some angular coordinate $\theta$. Then the metric:
\\
\\
\numberwithin{equation}{section}
{(20)}
$g=dr^2+r^{10}e^{(\lambda/2)r^2}d\theta^2$
\\
\\
The above case is singular and incomplete only at $r=0$ and its Gauss curvature is $K=-20/r^2-(11/2)\lambda-(\lambda^2/4)r^2$.
\\
\\
Also in this case for the initial hypothesis, $2K=R_{f_B}=cf$, we have $f=\frac{(11+\sqrt{57})\lambda}{8c}$ and $f=\frac{(11-\sqrt{57})\lambda}{8c}$. This means that $f$ is constant and this implies $h=0$. This case is not admitted in our analysis (we have set $h \neq 0$), so there is not solution for this case. Since if we consider a flat metric (i.e. $h=0$) this force $\lambda=0$ that is not the case in our interest. \; \; \; \; \; \; \;$\Box$
\\
\\
{\bfseries Remark} It is easy to verify that the results obtained in this paper are the same if we consider the Base-manifold with semi-Riemannian metrics.
\\
\\
{\bfseries Conclusions}: Under the examined conditions of $M$ and Fiber, the {\textit {$f$-curvature-Base}} is equivalent to requesting a flat metric on the Base-manifold. In fact for what we have seen, we can consider $M$ in only two types of Ricci-flat manifolds:
\begin{itemize}
\item [\ -] Einstein-Fiber (i.e. $\ddot{Ric}=\mu \ddot{g}$) with Flat-Base and $f$ affine,
\item [\ -] Flat-Fiber (i.e. $\ddot{Ric}=\mu \ddot{g}$, with $\mu=0$) with Flat-Base and $f$ constant. 
\end{itemize}
while cannot exist non-Ricci-flat $M$ with $R_{f_B}$ and Ricci-flat Fiber.
\\
\\
\\
{\bfseries References}
\\
\\
\; [1] Atiyah M. F. and Ward R. S. - "Instantons and algebraic geometry" - Commun. Math. Phys. 55(2) (1977), pp. 117-124.
\\
\; [2] Besse A. L. - "Einstein Manifolds" - Springer-Verlag (1987).
\\
\; [3] Castro C., Nieto J. A. - "On (2+2)-dimensional space time,strings and black holes" - International Journal of Modern Physics A 22(11) (2007) pp. 2021-2045.
\\
\; [4] De Andrade M. A., Del Cima O. M. and Colatto L. P. - "$N=1$ super-Chern-Simons coupled to parity-preserving matter from Atiyah-Ward space-time" - Phys. Lett. B 370(59) (1996), pp. 59-64.
\\
\; [5]  Evans L. C. - "Partial Differential Equations" Graduate Studies in Mathematics vol. 19 - American Mathematical Society (1997).
\\
\; [6] Ferrando J. J. and Saez J. A.- "An intrinsic characterization of 2+2 warped spacetimes" - Class. Quantum Grav. 27 (2010), 205023 (18 pages).
\\
\; [7] Kim D. S. - "Compact Einstein warped product spaces" - Information Center for Mathematical Sciences 5(2) (December 2002), pp. 1–5.
\\
\; [8] Kim D. S. and Kim Y. H. - "Compact Einstein warped product space with nonpositive scalar curvature" - Proceedings of the american mathematical society 131(8), pp. 2573-2576.
\\ 
\; [9]  Kuhnel W. Rademacherb H. - "Conformally Einstein product spaces" - Diff. Geometry and its App. 49 (December 2016), pp. 65-96.
\\
\; [10] Leandro B., Lemes de Sousa M. and Pina R. - "On the structure of Einstein warped product semi-Riemannian manifolds" -  Journal of Integrable Systems 3 (2018), pp. 1–11. 
\\
\; [11] Nieto J. A. - "Oriented matroid theory and loop quantum gravity in (2+2) and eight dimensions" - Rev. Mex. Fis. 57(2011) pp. 400–405.
\\
\; [12] Nieto J. A. and Pereyra C. - "Dirac equation in (1+3) and (2+2)-dimensions" - International Journal of Modern Physics A Vol. 28(24) (2013) 1350114 (18 pages).
\\
\; [13] O' Neill B. - "Semi-Riemannian geometry with applications to relativity" - Academic Press (1983). 
\\
\; [14] Serre D. - "Prolongement analytique et nombre de solutions d'une équation aux dérivées partielles elliptique non linéaire paramétrée" - C.R. Acad. Sci. Paris Sér. A (1978) 278, pp 1021-1023.
\\
\; [15] Ward R. S. - "Integrable and Solvable Systems, and Relations among them" - Phil. Trans. R. Soc. London, Ser. A 315 (1985), pp. 451-457.
\\

\par\bigskip
\end{document}